\documentclass[11pt,a4paper,twoside]{amsart}

\usepackage{amsfonts, indentfirst}
\usepackage{ifthen,amssymb}
\usepackage[all]{xypic}
\usepackage{young}
\usepackage{graphicx}
\usepackage{epsfig}
\usepackage[latin1]{inputenc}
\usepackage[T1]{fontenc}
\usepackage{mathptmx}
\usepackage{url}

\usepackage{tikz}
\usetikzlibrary{automata}

\newtheorem{prop}{Proposition}[section]
\newtheorem{theorem}[prop]{Theorem}

\newtheorem{defi}[prop]{Definition}
\newtheorem{coro}[prop]{Corollary}
\newtheorem{remark}[prop]{Remark}

\newcommand{\cqd}{\hfill$\Box$}

\title[AG codes, $t-$designs and partition sets]
{AG codes, $t-$designs and partition sets}

 \author[Cristina Mart{\'\i}nez]{Cristina Mart{\'\i}nez}
 \address{Hamilton Institute, 
Maynooth University, Maynooth, Co. Kildare, Ireland}
 \author[Alberto Besana]{Alberto Besana}

\address{FieldAware, Fumbally Square, Dublin 8, Ireland}
\email{cristina.martinezramirez@nuim.ie}
 
 \email{abesana@mat.uab.cat}

\begin{document}
\maketitle
\begin{abstract}
We study $q-$designs over finite fields and their connection to network coding. We prove that invariant subgrassmanians in the Grassmannian $\mathcal{G}_{k,n}(\mathbb{F}_{q})$ parametrizing $k-$planes in the vector space $\mathbb{F}_{q}^{n}$ by the action of any triangle group hold a $t$-design.

\end{abstract}

\section{Introduction}

Let $q$ be a power of a prime number $p$. It is well known, that there exists exactly one finite field with $q$ elements which is isomorphic to the splitting field of the polynomial $x^{q}-x$ over the prime field $\mathbb{F}_{p}$.

Any other field $F$ of characteristic $p$ contains a copy of $\mathbb{F}_{p}$. We denote respectively by $\mathbb{A}^{n}(\mathbb{F}_{q})$ and $\mathbb{P}^{n}(\mathbb{F}_{q})$ the affine space and the projective space over $\mathbb{F}_{q}$. Let $\mathbb{F}_{q}[x_{1},x_{2},\ldots, x_{n}]$ be the algebra of polynomials in $n$ variables over $\mathbb{F}_{q}$.

The encoding of an information word into a $k-$dimensional subspace is usually known as coding for errors and erasures in random network coding, \cite{KK}. 
      Namely, let $V$ be an $N-$dimensional vector space over $\mathbb{F}_{q}$, a code for an operator channel with ambient space $V$ is simply a nonempty collection of subspaces of $V$. The collection of subspaces is a code for error correcting errors that happen to send data through an operator channel.  The matrix coding the information is parametrised by random variables $a_{1}, a_{2},\ldots, a_{n}$ which constitute the letters of an alphabet. Here the operator channel is an abstraction of the operator encountered in random linear network coding, when neither transmitter nor receiver has knowledge of the channel transfer characteristics. The input and output alphabet for an operator channel is the projective geometry.       
 
Let $\mathcal{C}$ be a non-singular, projective, irreducible curve defined over  $\mathbb{F}_{q}$, as the vanishing locus of a polynomial $F\in \mathbb{F}_{q}[x_{0},x_{1},x_{2}]$.
We define the number $N(q)$ of $\mathbb{F}_{q}-$rational points on the curve to be
$$N(q)=|\{(x_{0},x_{1},x_{2})\in \mathbb{P}^{2}(\mathbb{F}_{q})| F(x_{0},x_{1},x_{2})=0\}|.$$
It is a polynomial in $q$ with integer coefficients, whenever $q$ is a prime power.

The number of points $\overline{\mathcal{C}}(\mathbb{F}_{q^{r}})$ on $\mathcal{C}$ over the extensions $\mathbb{F}_{q^{r}}$ of $\mathbb{F}_{q}$ is encoded in an exponential generating series, called the zeta function of $\overline{\mathcal{C}}$:
 $$Z(\mathcal{C},t)=exp\left( \sum_{r=1}^{\infty}\sharp \overline{\mathcal{C}}(\mathbb{F}_{q^{r}})\frac{t^{r}}{r}\right).$$
 
 One of the main problems in algebraic geometry codes is to obtain non-trivial lower bounds for the number $N(F_{i})$ of rational places of towers of function fields $\{F_{i}/\mathbb{F}_{q}\}_{i=1}^{\infty}$ such that $F_{i}\subsetneq F_{i+1}$.
 
For example, one can consider towers $\mathbb{F}_{q}(t^{\frac{1}{n}})$ or $\overline{\mathbb{F}}_{q}(t^{\frac{1}{n}})$, as $n$ varies through powers of the prime $p$ or through all integers not divisible by the characteristic of the ground field, that is $p$. Here $\overline{\mathbb{F}}_{q}$ denotes the algebraic closure of $\mathbb{F}_{q}$. The corresponding field extension $F=\mathbb{F}_{q}(t)\hookrightarrow F'=\mathbb{F}_{q}(t,y)$, where $y$ is a $n-$root of the polynomial $\sigma(u,t)=t^{n}-u \in F[u]$ and $u$ is transcendental over $\mathbb{F}_{q}$, is a finite cyclic Galois extension of degree $n$. Moreover $F'/F$ is an extension of Kummer type and $\mathbb{F}_{q}$ is the full constant field. 

Let $V$ be an $n+1$ dimensional vector space over the field $\mathbb{F}_{q}$, we denote by $PG(n,q)$ or $\mathbb{P}(V)$ the $n-$dimensional projective space over it. The set of all subspaces of dimension $r$ is the Grassmannian and it is denoted by $\mathcal{G}_{r,n}(\mathbb{F}_{q})$ or by $PG^{r}(n,q)$. Any code $C\subset \mathbb{F}_{q}^{n}$ of parameters $[n,k,d]$ is a point in the  Grassmannian $\mathcal{G}_{k,n}(\mathbb{F}_{q})$, where $d$ is the minimum distance under the Hamming Metric (HM).

 In general, for any integers $n, r$ with $n\geq r \geq 0$,  we call $\phi(r;n,q):=|PG^{r}(n,q)|$, the number of $r$ dimensional subspaces of an $n$ dimensional subspace over $\mathbb{F}_{q}$. It is  the number of ways of choosing $r+1$ linearly independent points in $PG(n,q)$ divided by the number of ways of choosing such a set of points in a particular $r-$space. It is given by the $q-$ary binomial coefficient,

 
 $$\phi(q;n,r)=\left[
\begin{matrix}
 n \\
 r\\
\end{matrix}
\right]_{q}=\frac{(q^{n+1}-1)(q^{n+1}-p)\ldots (q^{n+1}-q^{r})}{(q^{r+1}-1)(q^{r+1}-q)\ldots (q^{r+1}-q^{r})}.$$

\section{Algebraic geometric codes}
Let $\mathbb{F}_{q}$ be a finite field of $q$ elements, where $q$ is a power of a prime. The encoding of an information word into a $k-$dimensional subspace is usually known as coding for errors and erasures in random network coding \cite{KK}.       
 We consider as an alphabet a set $\mathcal{P}=\{P_{1},\ldots, P_{N}\}$ of $N-$ $\mathbb{F}_{q}$ rational points lying on a smooth projective curve $\mathcal{C}$ of genus $g$ and degree $d$  defined over the field $\mathbb{F}_{q}$. Goppa in \cite{Go}, constructed algebraic geometric linear codes from algebraic curves over finite fields with many rational places. 
      
\begin{defi} \label{AGC} Algebraic geometric codes (AGC) are constructed by evaluation of the global sections of a line bundle or more generally, a vector bundle on the curve $\mathcal{C}$ over $N$ $(N>g)$ distinct rational places $P_{1},\ldots, P_{N}$, being $g$ the genus of the curve. Namely, let $F| \mathbb{F}_{q}$ be the function field of the curve, 
 $\mathcal{D}$ the divisor $P_{1}+\cdots+ P_{N}$ and $G$ a divisor of $F| \mathbb{F}_{q}$ of degree $s\leq N$ such that $\rm{Supp}\,G\cap \rm{Supp}\, D=\emptyset$. Then the geometric Goppa code associated with the divisors $D$ and $G$ is defined by 
     $$\mathbf{ C}(D,G)=\{(x(P_{1}),\ldots, x(P_{N}))|\, x\in \mathcal{L}(G)\}\subseteq \mathbb{F}_{q^{n}}.$$
\end{defi}      
 
 Recall that $\mathbb{F}_{q^{n}}|\mathbb{F}_{q}$ is a cyclic Galois extension and 
is finitely generated by unique element $\alpha\in \mathbb{F}_{q}$. In particular, any element in $\mathbb{F}_{q^{n}}$ can be represented uniquely as a polynomial in $\alpha$ of degree less than $n$ with coefficients in $\mathbb{F}_{q}$. By Riemann-Roch Theorem, the code $\mathbf{ C}(D,G)$ of Definition \ref{AGC} is a linear $[n,k,d]$ code over $\mathbb{F}_{q}$, that is a code of length $n$, dimension $k$ with $k\geq s-g+1$, and minimum distance $d\geq n-s$. 
 
 Another important family of Goppa codes is obtained considering the rational normal curve $\mathcal{C}^{n}$ defined over $\mathbf{F}_{q}$:
 $$\mathcal{C}^{n}:=\{\mathbb{F}_{q}(1,\alpha, \ldots,\alpha^{n}):\, \alpha\in \mathbb{F}_{q}\cup \{\infty\}\}.$$
We assume that $q$ and $n$ are relatively prime $(n,q)=1$, so that any codeword $(c_{0},c_{1},\ldots, c_{n-1})$ can be expressed into a $q-$ary $k-$vector with respect to the basis $\{1,\alpha, \ldots, \alpha^{n}\}$. These codes are just Reed-Solomon codes or cyclic codes of parameters $[n,k,d]_{q}$ over $\mathbf{F}_{q}$ with parity check polynomial $h(x)=\prod_{i=1}^{q-1}(x-\alpha^{i})$ where $\alpha$ is a primitive root of $\mathbf{F}_{q}$ such that $\alpha^{k+1}=\alpha+1$. Recall that a linear cyclic code is an ideal in the ring $\mathbb{F}_{q}[x]/(x^{n}-1)$ generated by a polynomial $g(x)$ with roots in the splitting field $\mathbb{F}_{q}^{l}$ of $x^{n}-1$, where $n|\,q^{l}-1$, (\cite{BS2}). A natural question is how many polynomial are there over the algebraic closure of $\mathbb{F}_{q}[x]$. 
The next theorem expresses this number in terms of Stirling numbers.
 \begin{theorem} \label{number} The number of  polynomials decomposable into distinct linear factors over a finite field 
 $\mathbb{F}_{q^{n}}$ of arbitrary characteristic  a prime number $p$, is equal to
$\sum_{k=1}^{n}(q)_{k}$, where $(q)_{k}$ is the falling factorial polynomial  
$q\cdot (q-1)\ldots (q-(k-1))=\sum_{k=0}^{n}s(n,k)\,q^{k}$, where $s(n,k)$ is the Stirling number of the first kind 
divided by the order of the affine transformation group of the affine line
$\mathbb{A}^{1}=\mathbb{P}^{1}\backslash \infty$, that is $q^{2}-q$.
 
\end{theorem}
{\it Proof.} We need to count all the polynomials $f_{n}(x)$ in one variable of degree $n$ fixed. We assume that our polynomial $f_{n}(x)$ decomposes into linear factors, otherwise we work over $\bar{\mathbb{F}}_{q}[x]$, where $\overline{\mathbb{F}}_{q}$ denotes the algebraic closure of the finite field $\mathbb{F}_{q}$. Since the number of ordered sequences on $q$ symbols is $q!$ and each root is counted with its multiplicity, it follows that the number of monic polynomials with $n-1$ different roots is $q(q-1)(q-2)\ldots (q-n+1):=(q-2)_{n}$. Now we observe that polynomials are invariant by the action of automorphisms of the affine line, so we must divide this number by the order of this group which is $q^{2}-q$.
\cqd

\begin{theorem} \label{generators} Given a set of integers $\{0,1,\ldots, n-1\}$ modulo $n$, there is a set $J$  of $k$ integers which is a set of roots, that is, there is a polynomial $h(x)=\prod_{j\in J}(x-\alpha^{j})$, where $\alpha$ is a generator of $(\mathbb{F}_{p})^{m}$ for some prime number $p$ and $m$ is the least integer such that $n|p^{m}-1$. The ideal $h(x)$ generates in $\mathbb{F}_{p^{m}}[x]/(x^{n}-1)$ is a cyclic linear code of parameters $(n,k, n-k+1)$.
\end{theorem}

\begin{remark}
As an application of Theorem \ref{number} and Theorem \ref{generators}, given an integer $n$, we can count the number of cyclic codes of parameters $[n,k]$ for each $0\leq k \leq n$ and set of roots $\alpha_{1},\ldots, \alpha_{k}$ in the splitting field of $x^{n}-1$, the corresponding polynomial $g(x)=\prod_{i=1}^{k}(x-\alpha_{i})$ generates a linear cyclic code in the ring $\mathbb{F}_{q}[x]/(x^{n}-1)$. Thus for each $0\leq k\leq n$ there are exactly $(q)_{k}/(q^{2}-q)$ cyclic codes.

\end{remark}

\section{$t-$designs and representation theory of $GL(n,\mathbb{F}_{q})$}
A simple $t-$design over a finite field or, more precisely, a $t-(n,k,\lambda;q)$ design is a set $\mathcal{B}$ of $k-$subspaces of an $n-$dimensional vector space $V$ over the finite field $\mathbb{F}_{q}$ such that each $t-$subspace of $V$ is contained in exactly $\lambda$ blocks of $\mathcal{B}$.
Recently, designs over finite fields gained a lot of attention because of its applications for error-correcting in networks. 
If $V$ is the finite field $\mathbb{F}_{q}^{n}$, then the set of points are the vectors and the block set $\mathcal{B}$ of $k$ subspaces $K\subseteq \mathbb{F}_{q}^{n}$ are the points in the Grassmannian  $\mathcal{G}_{k,n}(\mathbb{F}_{q})$.
A permutation matrix $\sigma \in GL(n,q)$ acts on the Grassmannian by multiplication on the right of the corresponding representation matrix. In particular $\sigma$ is an automorphism of the design $\mathcal{D}=(V,\mathcal{B})$ if and only if $\sigma$ leaves the Grassmannian invariant, that is $\mathcal{B}^{\sigma}=\mathcal{B}$. 

In particular, we are interested in understanding the orbits by the action of any permutation matrix of $GL(n,q)$ and moreover of any subgroup $G$ contained in $GL(n,q)$. Further, it is possible to count the orbits of the action in several cases and these correspond to blocks of the design satisfying certain geometrical properties. 


\begin{defi} Let $\alpha\in \mathbb{F}_{q^{n}}$ be a generator of the underlying vector space over $\mathbb{F}_{q}$. Then an $r-$dimensional $W$  subspace is $\alpha-$splitting if $\alpha^{i}W=W$ is invariant under the action of any element $\alpha^{i}$ in the Galois group of the extension $\mathbb{F}_{q}\hookrightarrow \mathbb{F}_{q}(\alpha)$. 

More precisely, given any $\mathbb{F}_{q}-$linear endomorphism $T:\,\mathbb{F}_{q^{n}}\rightarrow \mathbb{F}_{q^{n}}$, an $r-$di\-men\-sio\-nal subspace is $T-$splitting if $\mathbb{F}_{q^{n}}=W\oplus T(W)\oplus \cdots \oplus T^{n-1}(W),$ where $T^{j}$ denotes the $j-$fold composite of $T$ with itself.
\end{defi}

\begin{prop} Let $\lambda$ be the number of $\alpha-$splitting subspaces of $\mathbb{F}_{q}^{n}$, then 
the design whose blocks are  the $\alpha-$splitting subspaces of $\mathbb{F}_{q}^{n}$  is a $t-(n,r,\lambda;q)$ design.
\end{prop}

{\it Proof.} Suppose $n=r\cdot s$, where $r, s$ are coprime positive integer numbers, that is $(r,s)=1$. 
There is an element $\alpha\in \mathbb{F}_{q}^{n}$ of order $s$, so that $\{1, \alpha, \ldots, \alpha^{s-1}\}$ is a basis of $\mathbb{F}_{q^{s}}$ over $\mathbb{F}_{q}$. Hence the set $\{1, \alpha^{s}, \alpha^{2s}\ldots, \alpha^{(r-1)s}\}$ spans an $r-$dimensional $\alpha-$splitting subspace of $\mathbb{F}_{q^{n}}$, say $W$. 
Define an isomorphism between $\mathbb{F}^{rs}_{q}$ and $\mathbb{F}^{r}_{q^{s}},$ for $i=0,1,\ldots, ,r-1$ associating each $s-$tuple $(v_{i,0},v_{i,1},\ldots, v_{i,s-1}),$ with the element $v_{i}\in \mathbb{F}_{q^{s}},$ where $v_{i}=v_{i,0}+v_{i,1}\alpha+\ldots +v_{i,s-1}\alpha^{s-1}$. Then every element in $\mathbb{F}^{rs}_{q}$ is in correspondence with an element in $\mathbb{F}^{r}_{q^{s}}$.
If we complete $W$ to a basis of $\mathbb{F}_{q^{n}}$ by adding elements $v_{1},\ldots, v_{m}\in \mathbb{F}_{q^{n}}$, we let 
$$\mathcal{B}^{\alpha}_{(v_{1},\ldots, v_{m})}:=\{v_{1},\ldots, v_{m},\alpha v_{1},\ldots, \alpha v_{m},\ldots, \alpha^{n-1}v_{1},\ldots, \alpha^{n-1}v_{m}\},$$ where $\mathcal{B}^{\alpha}_{(v_{1},\ldots, v_{m})}$, is regarded as an ordered set with $n=r\cdot m$ elements.  In this case $\mathcal{B}^{\alpha}_{(v_{1},\ldots, v_{s})}$ is necessarily a $\mathbb{F}_{q}-$basis of an $s-$dimensional subspace of $\mathbb{F}_{q^{rs}},$ and we will refer to $\mathcal{B}^{\alpha}_{(v_{1},\ldots, v_{s})}$ as an $\alpha-$splitting ordered basis of $\mathbb{F}_{q}^{rs}$. The subspaces generated by the $\alpha-$splitting ordered basis constitute an $s-$design of parameters $(n,k,\lambda;q)$, where $k$ is the number of ordered basis, which is exactly $|{\rm{GL}}(m,\mathbb{F}_{q})|=\prod_{i=0}^{m-1}(q^{m}-q^{i})$.

\cqd

Let us define $S(\alpha, r,s;q)$ to be the number of $\alpha-$splitting subspaces of $\mathbb{F}_{q^{n}}$ of dimension $r$. 
Let $N(\alpha, r, s;q)$ be the number of ordered basis of $\mathbb{F}_{q^{rs}}$, then 
$$S(\alpha, r,s;q)=\frac{N(\alpha, r,s;q)}{|{\rm{GL}}_{n}(m,\mathbb{F}_{q})|}.$$
\begin{coro}  If $\mathbb{F}_{q^{n}}=\mathbb{F}_{q}(\alpha)$, any Galois conjugate $\beta$ of $\alpha$ generates $\mathbb{F}_{q^{n}}$ and the corresponding $t-$designs are isomorphic.
\end{coro}

For example, if $\beta=\frac{a\alpha^{q^{r}}+b}{c\alpha^{q^{r}}+d}$ for some non negative integer $r$ and 
$\begin{pmatrix}  a  &  b  \\  c & d  \end{pmatrix}\in GL_{2}(\mathbb{F}_{q})$, then $\mathbb{F}_{q^{n}}=\mathbb{F}_{q}(\beta)$ and $S(\alpha, r,s;q)=S(\beta,r,s;q)$.


We want to understand which subspaces are invariant by the action of elements of the general linear group $GL(n, \mathbb{F}_{q})$ or finite subgroups of $GL(n,\mathbb{F}_{q})$. In this way, one can construct designs with prescribed groups where the blocks are the orbits by the action, and thus to generalize to other Galois extensions not necessarily cyclic.

The general linear group $GL(n,\mathbb{F}_{q})$ acts on the Grassmannian $\mathcal{G}_{k,n}(\mathbb{F}_{q})$ by multiplication on the right:
 \begin{eqnarray}
 \mathcal{G}_{k,n}(\mathbb{F}_{q})\times GL(n,\mathbb{F}_{q}) \rightarrow \mathcal{G}_{k,n}(\mathbb{F}_{q})\\
 (\mathcal{U},A)\rightarrow \mathcal{U}A,
 \end{eqnarray} for a representation matrix $\mathcal{U}\in \mathbb{F}_{q^{k\times n}}.$
 
 Then cyclic codes correspond to orbits in $\mathcal{G}_{k,n}(\mathbb{F}_{q})$ defined by a cyclic subgroup $\mathbb{Z}_{p}$ for some prime number $p$. These codes are supported on the Normal Rational Curve (NRC), that is, coding vectors for networks with $n$ sources live in the projective space $PG(n-1,q)$.
 

\begin{defi}
In $PG(n-1,q)$ a $(k;r)-$arc is a set of $k$ points any $r$ of which form a basis for $\mathbb{F}_{q}^{n}$, or in other words, $r-1$ of them but not $r$ are collinear.
\end{defi}

Consider the normal rational curve over $\mathbb{F}_{q}$:
$$\mathcal{V}^{n}_{1}:=\Big\{\mathbb{F}_{q}(1,x,x^{2},\ldots,x^{n})| \ x\in \mathbb{F}_{q}\bigcup \{\infty\}\Big\}$$ is a $(q+1)-$arc in the $n-$dimensional projective space $PG(n,q)$.

In coordinates, the Veronese mapping

$\mathbb{F}_{q}(x_{0}b_{0}+x_{1}b_{1})\rightarrow \mathbb{F}_{q}(\sum_{l=0}^{n}x_{0}^{n-e}x_{1}^{e}c_{e}) \ \ x_{i}\in \mathbb{F}_{q}$, maps the point set of the projective line $\mathcal{P}(X)$ into the point set of $\mathcal{P}(Y)$,  where $Y$ is an $(n+1)-$dimensional vector space over $\mathbb{F}_{q}$ with a basis $\{c_{0},\ldots, c_{n}\}$ where $n\geq 2$.

We see that if $q\leq n$, the NRC is a basis of a $q-$dimensional projective subspace, that is, a $PG^{q}(n,q)$. So we can enumerate how many NRC's  are there in a $PG(n,q)$. The answer is $\phi(q;n,q)$, the number of ways of choosing such a set of points in a particular $q-$space.
If $q\geq n+2$, the NRC is an example of a $(q+1)-$arc. It contains $q+1$ points, and every set of $n+1$ points are linearly independent.

When $q=2$, finite subgroups of $PGL(2,q)$ have been classified. They are isomorphic either to the dihedral group $D_{n}$ of order $2n \,(n\geq 2)$, the alternating group $A_{4}$, the alternating group $A_{5}$, or the symmetric group $S_{4}$. The corresponding invariant subgrassmannianas in $\mathcal{G}_{k,n}(\mathbb{F}_{2}^{m})$ define $t-$designs.

\begin{prop}\label{PropSteiner}
Invariant subgrassmannians in $\mathcal{G}_{k,n}(\mathbb{F}_{2}^{m})$ by the action of any triangle group hold a $t-$design.
\end{prop}
{\it Proof.} Triangle groups are reflection groups which admit a presentation
$${\rm{T}}(r,s,t)=\langle x, y, x:\, x^{r}=y^{m}=z^{w}=xyz=1\rangle,$$ with $r, m, w$ integer numbers such that $\frac{1}{r}+\frac{1}{m}+\frac{1}{w}<1$.
They are finite subgroups of ${\rm{PGL}}(2,q)$ and they are known to be isomorphic to either the dihedral group $D_{m}$ of order $2m, m\geq 2$, the alternating group $A_{4}$, the alternating group $A_{5}$, or the symmetric group $S_{4}$. Since they are finitely generated, the invariant subgrassmannians in $\mathcal{G}_{k,n}(\mathbb{F}_{2}^{m})$ define a $t-$design where $t$ is the number of generators. The dihedral group $D_{m}$ is generated by a rotation $\tau$ and a reflection $\sigma$, then the corresponding invariant subgrassmannian in $\mathcal{G}_{k,n}(\mathbb{F}_{2}^{m})$ define a 2-design and so on for the other triangle groups.

\cqd

\begin{remark} From a graph theoretical point of view, we associate to the 2-design generated by $\tau$ and $\sigma$ in Proposition \ref{PropSteiner} the graph which has as vertex set $V$ the points of the projective system $\mathbb{P}((\mathbb{F}_{2})^{m})$ and edge set $E\subseteq [V]^{2}$ the lines of the projective space which corresponds to the blocks of the design. There are $\left[
\begin{matrix}
 m \\
 2\\
\end{matrix}
\right]_{q}$ lines. For any two points there are as much blocks (lines) containing these points as eigenvectors $W_{j}$  by the action of the linear operators $\tau$ and $\sigma$. This special design with parameters $t=2$ and $k=3$ is a Steiner triple system.
\end{remark}
\section{Conclusion}
The proofs given in this short note give a way of constructing designs of given parameters $n, q,  s$ by deriving ordered bases of $(\mathbb{F}_{q})^{n}$.

 \end{document}